\documentclass[a4paper,10pt]{amsart}
\usepackage[utf8]{inputenc}

\usepackage{amsmath, amssymb, amsthm}
\theoremstyle{remark}
\newtheorem{ex}{Example}

\title[A Macaulay2 package for characteristic classes]{A Macaulay2 package for characteristic classes and the topological Euler characteristic \\of complex projective schemes}
\author{Christine Jost}

    \address{Stockholm University, Department of
             Mathematics, SE-106 91 Stockholm, Sweden}
    \email{jost@math.su.se}
    \urladdr{http://www.math.su.se/$\sim$jost}

\thanks{}

\begin{document}

\begin{abstract}
The Macaulay2 package \emph{CharacteristicClasses} provides commands for the computation of the topological Euler characteristic, the degrees of the Chern classes and the degrees of the Segre classes of a closed subscheme of complex projective space. The computations can be done both symbolically and numerically, the latter using an interface to Bertini. We provide some background of the implementation and show how to use the package with the help of examples. 
\end{abstract}

\maketitle

\section{Introduction}

%
The Macaulay2 \cite{M2} package \emph{CharacteristicClasses} computes degrees of Chern and Segre classes of complex projective schemes. It also provides a command computing the topological Euler characteristic.
Recall that the Chern classes of a smooth scheme are defined to be the Chern classes of the tangent bundle. Moreover, the Segre classes of a scheme $X$ embedded in $\mathbb{P}^n$ are defined to be the Segre classes of the normal cone $C_X\mathbb{P}^n$ of the embedding.  Note that the degree of a $d$-dimensional cycle class, i.e., the weighted sum of classes of $d$-dimensional subvarieties, is defined to be the weighted sum of the degrees of the subvarieties.

The computations done by the package \emph{CharacteristicClasses} are based on a number of related algorithms which reduce the problem to the computation of the degrees of residuals, i.e., components complementary to a certain scheme in the intersection of hypersurfaces containing that scheme. The algorithm computing Chern classes is described in \cite{DEPS}, the one computing Segre classes for possibly singular schemes in \cite{EJP}, and the computation of the topological Euler characteristic is described in \cite{J}. The residuals can be computed either symbolically or nu\-mer\-ically, which yields a symbolical and a numerical version of the algorithms. Both versions are implemented in the package, the numerical version via an interface to Bertini \cite{B}. 

There are other ways of computing the topological Euler characteristic and degrees of characteristic classes. The Macaulay2 command \emph{euler} computes the topological Euler characteristic of smooth projective varieties by computing the Hodge numbers. In \cite{A}, Aluffi describes algorithms for the computation of degrees of characteristic classes and the topological Euler characteristic. These algorithms are implemented in Aluffi's package \emph{CSM}, which is not part of Macaulay2 but can be obtained at \texttt{www.math.fsu.edu/$\sim$aluffi/CSM/CSM.html}. The different implementations complement each other, as is shown in \cite{EJP}.

This article describes version 0.2 of the package \emph{CharacteristicClasses}, which can be obtained at \texttt{www.math.su.se/$\sim$jost/CharacteristicClasses.m2}. Earlier versions provide computations of degrees of Chern and Segre classes, but not of the topological Euler characteristic.

\section{Computing degrees of Chern and Segre classes}

We describe how to make Chern and Segre classes of projective schemes computationally tractable by computing their degrees. The standard reference for the notions used in this section, Chern classes, Segre classes and Chow groups, is \cite{F}.

Let $X$ be a $k$-dimensional closed subscheme of $\mathbb{P}^n$, embedded by $i\colon X \hookrightarrow \mathbb{P}^n$, and denote its Chow group by $A_*(X) = \bigoplus_{d=0}^kA_d(X)$. Characteristic classes of $X$ are elements of the graded group $A_*(X)$, whose generators are generally hard to compute. Hence algorithms for the computation of characteristic classes focus on computing a coarser invariant, the degrees of the classes. Let $\alpha = \sum_i \alpha_i [V_i]$ be a cycle class in $A_d(X)$, given as the weighted sum of classes of $d$-dimensional subvarieties $V_i$. Then the degree $\mathrm{deg}(\alpha)$ of $\alpha$ is defined to be the weighted sum $\sum_i \alpha_i \, \mathrm{deg}(V_i)$ of the degrees of the varieties $V_i$, seen as subvarieties of $\mathbb{P}^n$. One can also consider the pushforward of the cycle $\alpha$ to the Chow group of  $\mathbb{P}^n$, which is known to be $A_*(\mathbb{P}^n) = \mathbb{Z}[H]/(H^{n+1})$, where $H$ is the class of a general hyperplane. The two viewpoints are 
equivalent because 
$i_*(\alpha) = \mathrm{deg}(\alpha)H^{n-d}$. Both points of view, 
degrees of the classes and pushforward to the Chow ring of $\mathbb{P}^n$, are supported by the package \emph{CharacteristicClasses}.

The Chern classes $c_1(X), \ldots, c_k(X)$ of a smooth scheme $X$ are by definition the Chern classes $c_1(T_X), \ldots, c_k(T_X)$ of the tangent bundle of $X$. The total Chern class is the sum $1 + c_1(X) + \ldots + c_k(X)$. Furthermore, the Segre classes $s_1(X,\mathbb{P}^n), \ldots,$ $s_k(X,\mathbb{P}^n)$ of a possibly singular $X$ are the Segre classes $s_1(C_X\mathbb{P}^n), \ldots, s_k(C_X\mathbb{P}^n)$ of the normal cone $C_X\mathbb{P}^n$ of $X$ in $\mathbb{P}^n$. If the embedding is regular, the normal cone is a vector bundle, called the normal bundle. The total Segre class is the sum $1 + s_1(X,\mathbb{P}^n) + \ldots + s_k(X,\mathbb{P}^n)$.

The package \emph{CharacteristicClasses} provides the commands \emph{chernClass} and \emph{segre\-Class} which compute the degrees of the Chern and Segre classes of a given closed subscheme of $\mathbb{P}^n$ or, equivalently, the pushforward of the total Chern and Segre class to the Chow group of $\mathbb{P}^n$. The closed subscheme $X$ is given by a number of generators of a homogeneous ideal in a polynomial ring. 
The main idea of the algorithm is to relate the degrees of the Chern and Segre classes to the degrees of so-called residuals. One chooses randomly a number of hypersurfaces containing the scheme $X$, where the number of hypersurfaces is at least the codimension of $X$. According to a Bertini-type theorem, they intersect in $X$ and some possibly empty component of expected codimension, the residual. Its degree can be computed either symbolically or numerically. Symbolically one computes  the saturation using 
Gröbner basis techniques. The numerical computations can be done using software for the numerical solution of polynomial equation systems. The package \emph{CharacteristicClasses} implements both the symbolic version of the algorithms and the numerical, using an interface to Bertini \cite{B} for the latter. More details on the algorithms can be found in \cite{DEPS} and \cite{EJP}. 

\section{Computing the topological Euler characteristic}

The topological Euler characteristic of a complex projective variety is the Euler characteristic of the underlying topological space with the usual Euclidean topology, i.e., the alternating sum of the Betti numbers. For the computations with the package \emph{CharacteristicClasses}, we use that the topological Euler characteristic is equal to the degree of the top Chern-Schwartz-MacPherson class. More generally, the package can also compute the degrees of all the Chern-Schwartz-MacPherson classes. It uses the algorithm described in \cite{J}, which reduces the computation of the degrees of Chern-Schwartz-MacPherson classes to the computation of Segre classes of singular subvarieties. For a more detailed introduction to Chern-Schwartz-MacPherson classes we refer to \cite{A} and \cite{J}.

\section{Using the package CharacteristicClasses}

We present three examples which demonstrate how to use the package \emph{CharacteristicClasses}. Observe that this article describes version 0.2 of the package \emph{CharacteristicClasses}, but only version 0.1 is included in a standard installation of version 1.5 of Macaulay2, the most recent version at the moment of writing this article. More detailed information on the commands provided by this package can be found in its documentation, which is shown by e.g.\ using the command \emph{viewHelp CharacteristicClasses}.

\begin{ex}

We compute the degrees of the Chern classes of a toy example, the twisted cubic. We start by loading the package and defining the ideal of the twisted cubic, which is generated by the 2-by-2 minors of a 2-by-3 matrix. Then we use the command \emph{chernClassList} to obtain a list of the degrees of the Chern classes of the twisted cubic. The twisted cubic $C_\mathrm{tw}$ has dimension 1, hence it only has one Chern class $c_1(C_\mathrm{tw})$. The output of the command \emph{chernClassList} is the list $\{\mathrm{deg}(C_\mathrm{tw}), \mathrm{deg}(c_1(C_\mathrm{tw})) \}$.

\begin{verbatim}
Macaulay2, version 1.5
with packages: ConwayPolynomials, Elimination, IntegralClosure, 
               LLLBases, PrimaryDecomposition, ReesAlgebra, 
               TangentCone

i1 : loadPackage "CharacteristicClasses";
--loading configuration for package "CharacteristicClasses" from 
      file /.../.Macaulay2/init-CharacteristicClasses.m2

i2 : R = QQ[x,y,z,w];

i3 : twistedCubic = minors(2,matrix{{x,y,z},{y,z,w}})

               2                        2
o3 = ideal (- y  + x*z, - y*z + x*w, - z  + y*w)

o3 : Ideal of R

i4 : chernClassList twistedCubic

o4 = {3, 2}

o4 : List
\end{verbatim}

The computations above confirm that the twisted cubic has degree 3. As the twisted cubic is smooth, the degree of the top Chern class equals its Euler characteristic $\chi(C_\mathrm{tw})$, which is related to the genus $g(C_\mathrm{tw})$ by $\chi = 2-2g$. As $\mathrm{deg}(c_1(C_\mathrm{tw}))=2$ by the computations above, they confirm that the genus is 0 and the twisted cubic is a rational curve. 

As said in Section 2, computing the degrees of the Chern classes is equivalent to computing the push-forward of the total Chern class to the Chow ring of $\mathbb{P}^n$. The total Chern class of the twisted cubic is $1 + c_1(C_\mathrm{tw})$, its pushforward to the Chow ring of projective space is $\mathrm{deg}(C_\mathrm{tw}) H^2 + \mathrm{deg}(c_1(C_\mathrm{tw}))H^3$, where $H$ is the hyperplane class. The pushforward can be computed using the command \emph{chernClass}.
\newpage

\begin{verbatim}
i5 : chernClass twistedCubic 

       3     2
o5 = 2H  + 3H

     ZZ[H]
o5 : -----
        4
       H
\end{verbatim}

All computations can also be done numerically using an interface to Bertini. Version 1.3 of Bertini needs to be installed and the package must be configured correctly. For more information on the configuration, use \emph{viewHelp "configuring Bertini"}. Then
by using the value \emph{Bertini} for the option \emph{ResidualStrategy}, the computations are done numerically using Bertini instead of Gröbner basis computations.

\begin{verbatim}
i6 : chernClassList(twistedCubic, ResidualStrategy=>Bertini)

o6 = {3, 2}

o6 : List
\end{verbatim}
\end{ex}

\begin{ex}
We proceed with an example for the computation of Segre classes. The Whitney umbrella is a singular surface in $\mathbb{P}^3$. Over the reals its looks like an umbrella and the singular locus like the handle of the umbrella. We compute the degree of the first Segre class $s_1(S,\mathbb{P}^n)$ of the singular locus $S$. Observe that in the following the singular locus is represented in Macaulay2 not as an ideal, but as a projective variety. All commands in this package work for both ideals and projective varieties.

\begin{verbatim}
i7 : whitney = ideal(x^2*w - y^2*z);

o7 : Ideal of R

i8 : handle = Proj singularLocus whitney;

i9 : segreClassList handle

o9 = {1, 0}

o9 : List
\end{verbatim}

We see that the singular locus is a curve with degree 1, hence a line, and the degree of its first Segre class is $\mathrm{deg}(s_1(S,\mathbb{P}^n)) = 0$. According to \cite{A3}, the degree of the first Segre class is $\mathrm{deg}(s_1(S,\mathbb{P}^n)) = \nu - 2$, where $\nu$ is the number of pinch points of the surface. It follows that the Whitney umbrella has two pinch points. Only one of them is visible in the real patch $\{w\neq0\}$ in which the Whitney umbrella usually is drawn, it is the tip of the umbrella.
\end{ex}

\begin{ex}
We continue with an example from algebraic statistics which uses computations of the topological Euler characteristic, Example 2.2.2 from \cite{DSS}. The theory behind the computations is described in more detail in \cite{J}. The random censoring model with two events is a statistical model implicitely described by the ideal $(2p_0p_1p_2 + p_1^2p_2 + p_1p_2^2 - p_0^2p_{12} + p_1p_2p_{12})$ in the polynomial ring $\mathbb{C}[p_0,p_1,p_2,p_{12}]$, where $p_0$, $p_1$, $p_2$ and $p_{12}$ describe the probabilities of two events to occur before or after a third event. Given experimental data one would like to compute the values for the probabilities $p_0$, $p_1$, $p_2$ and $p_{12}$ which describe the data best, by maximizing the likelihood function. However, the function may have several stationary points and methods such as the Newton method may only find a local maximum. Hence it makes sense to define the maximum likelihood degree as the number of critical points of the likelihood function, 
as was done in \cite{CHKS}. By a 
theorem of Huh \cite{H}, for a large class of examples the maximum likelihood degree equals the signed topological Euler characteristic of a certain open subvariety of the model. In this case, it is the topological Euler characteristic of the open subvariety $V(2p_0p_1p_2 + p_1^2p_2 + p_1p_2^2 - 
p_0^2p_{12} + p_1p_2p_{12}) \setminus V(p_0p_1p_2p_{12}(p_0+p_1+p_2+p_{12}))$, where $p_0p_1p_2p_{12}\neq 0$ means that no probability should be zero, and $p_0+p_1+p_2+p_{12} \neq 0$ means that the probabilities should sum up to 1.  We compute the topological Euler characteristic using the command \emph{eulerChar} together with the inclusion-exclusion principle.

\begin{verbatim}
i10 : S = QQ[p0,p1,p2,p12];

i11 : randomCensoring = ideal(2*p0*p1*p2 + p1^2*p2 + p1*p2^2 - 
	      p0^2*p12 + p1*p2*p12);

o11 : Ideal of S

i12 : boundary = ideal( p0*p1*p2*p12*(p0+p1+p2+p12) ) + 
	      randomCensoring;

o12 : Ideal of S

i13 : eulerChar randomCensoring

o13 = 5

i14 : eulerChar boundary

o14 = 2
\end{verbatim}

It follows that the topological Euler characteristic of $V(2p_0p_1p_2 + p_1^2p_2 + p_1p_2^2 - 
p_0^2p_{12} + p_1p_2p_{12}) \setminus V(p_0p_1p_2p_{12}(p_0+p_1+p_2+p_{12}))$ is $5-2=3$. Hence the maximum likelihood degree of the random censoring model is 3, which confirms the result in \cite{DSS}.
 
\end{ex}

\section*{Acknowledgements}
Many thanks to the organizers and participants of the Macaulay2 workshop 2011 at the IMA in Minneapolis, where parts of the package were written. Also many thanks to my advisor Sandra Di Rocco for help with the article.

Jon Hauenstein adapted the output of the regenerative cascade in Bertini to the interface used in \emph{CharacteristicClasses}. 
The example for the computation of the number of pinch points of the Whitney umbrella is due to David Eklund.


\begin{thebibliography}{999}

\bibitem{A3} P. Aluffi, \emph{MacPherson's and Fulton's Chern classes of hypersurfaces}, International Mathematics Research Notices 11 (1994), 455--465.

\bibitem{A} P. Aluffi, \emph{Computing characteristic classes of
  projective schemes}, Journal of Symbolic Computation 35 (2003),
  3--19.

\bibitem{B} D. J. Bates, J. D. Hauenstein, A. J. Sommese, C. W. Wampler, \emph{Bertini: Software for Numerical Algebraic Geometry}, available at http://www.nd.edu/$\sim$sommese/bertini.

\bibitem{CHKS} F. Catanese, S. Hosten, A. Khetan, B. Sturmfels, \emph{The maximum likelihood degree}, American Journal of Mathematics 128 (2006), 671--697. 


\bibitem{DEPS} S. Di Rocco, D. Eklund, C. Peterson, A.J. Sommese, \emph{Chern numbers of smooth varieties via homotopy continuation and intersection theory}, Journal of Symbolic Computation 46 (2011), 23--33.

\bibitem{DSS} M. Drton, B. Sturmfels and S. Sullivant, \emph{Lectures on algebraic statistics}, Birkh{\"a}user Basel, 2008.

\bibitem{EJP} D. Eklund, C. Jost, C. Peterson, \emph{A method to compute Segre classes of subschemes of projective space}, Journal of Algebra and its Applications 12(2) (2013)

\bibitem{F} W. Fulton, \emph{Intersection theory}, Springer, 1984.

\bibitem{M2} D. R. Grayson, M. E. Stillman, \emph{Macaulay2, a software system for research in algebraic geometry}, available at http://www.math.uiuc.edu/Macaulay2/.


\bibitem{H} J. Huh, \emph{The maximum likelihood degree of a very affine variety}, arXiv:1207.0553 [math.AG].

\bibitem{J} C. Jost, \emph{An algorithm for computing the topological Euler characteristic of complex projective varieties}, arXiv:1301.4128 [math.AG].


\end{thebibliography}
\end{document}